\documentclass[11pt]{article}
\usepackage{euscript,amsmath,amssymb }
\usepackage{srcltx}

\title{On transverse hyperplanes to self-similar Jordan arcs.}

\author{Andrey Tetenov\footnote{Supported by Russian Foundation of Basic Research project
13-01-00513.}}

\begin{document}

\newcommand{\rr}{\mathbb{R}}
\newcommand \nn {\mathbb{N}}
\newcommand \zz {\mathbb{Z}}
\newcommand \bbc {\mathbb{C}}
\newcommand \rd {\mathbb{R}^d}

 \newcommand {\al} {\alpha}
\newcommand {\be} {\beta}
\newcommand {\da} {\delta}
\newcommand {\Da} {\Delta}
\newcommand {\ga} {\gamma}
\newcommand {\Ga} {\Gamma}
\newcommand {\la} {\lambda}
\newcommand {\La} {\Lambda}
\newcommand{\om}{\omega}
\newcommand{\Om}{\Omega}
\newcommand {\sa} {\sigma}
\newcommand {\Sa} {\Sigma}
\newcommand {\te} {\theta}
\newcommand {\fy} {\varphi}
\newcommand {\ep} {\varepsilon}
\newcommand{\e}{\varepsilon}

\newcommand{\VEC}{\overrightarrow}
\newcommand{\IN}{{\subset}}
\newcommand{\NI}{{\supset}}
\newcommand \dd  {\partial}
\newcommand {\mmm}{{\setminus}}
\newcommand{\probel}{\vspace{.5cm}}
\newcommand{\8}{{\infty}}
\newcommand{\0}{{\varnothing}}
\newcommand{\vse}{$\blacksquare$}

\newcommand {\bfep} {{{\bar \varepsilon}}}
\newcommand {\Dl} {\Delta}
\newcommand{\vA}{{\vec {A}}}
\newcommand{\vB}{{\vec {B}}}
\newcommand{\vF}{{\vec {F}}}
\newcommand{\vf}{{\vec {f}}}
\newcommand{\vh}{{\vec {h}}}
\newcommand{\vJ}{{\vec {J}}}
\newcommand{\vK}{{\vec {K}}}
\newcommand{\vP}{{\vec {P}}}
\newcommand{\vX}{{\vec {X}}}
\newcommand{\vY}{{\vec {Y}}}
\newcommand{\vZ}{{\vec {Z}}}
\newcommand{\vx}{{\vec {x}}}
\newcommand{\va}{{\vec {a}}}
\newcommand{\vga}{{\vec {\gamma}}}

\newcommand{\eS}{{\EuScript S}}
\newcommand{\eH}{{\EuScript H}}
\newcommand{\eC}{{\EuScript C}}
\newcommand{\eP}{{\EuScript P}}
\newcommand{\eT}{{\EuScript T}}
\newcommand{\eG}{{\EuScript G}}
\newcommand{\eK}{{\EuScript K}}
\newcommand{\eF}{{\EuScript F}}
\newcommand{\eZ}{{\EuScript Z}}
\newcommand{\eL}{{\EuScript L}}
\newcommand{\eD}{{\EuScript D}}
\newcommand{\E}{{\EuScript E}}
\def \diam {\mathop{\rm diam}\nolimits}
\def \fix {\mathop{\rm fix}\nolimits}
\def \Lip {\mathop{\rm Lip}\nolimits}

\newtheorem{thm}{\bf Theorem}
 \newtheorem{cor}[thm]{\bf Corollary}
 \newtheorem{lem}[thm]{\bf Lemma}
 \newtheorem{prop}[thm]{\bf Proposition}
 \newtheorem{dfn}[thm]{\bf Definition}

\newcommand{\dok}{{\bf{Proof}}}

\maketitle

\begin{abstract} We consider self-similar Jordan arcs $\ga$ in $\rd$, different from a line segment
 and  show that  they cannot be projected to a line bijectively. Moreover, we show that  the
 set of points $x\in\ga$, for which there is a hyperplane,
 intersecting $\ga$ at the point $x$ only, is nowhere  dense in $\ga$.
\end{abstract}
\maketitle
\noindent MSC classification: Primary 28A80\vspace{12mm}\\
Andrey Tetenov\\
Gorno-Altaisk state university\\
 649000 Gorno-Altaisk, Russia\\
e-mail:  {atet@mail.ru}

\section{Introduction.}
The first examples of  self-similar fractals  which appeared in
the beginning of XX century were the constructions of self-similar
curves with predefined geometrical properties \cite{Koch, Serp}.
Though   the study of geometrical properties of self-similar
curves is so close to historical origins of fractal geometry, some
of their elementary geometric properties were established only in
recent times.

For example, it was a common opinion   that self-similar curves
have no tangent at any of their points. But in 2005 A.Kravchenko
\cite{Krav} found that there are self-affine curves which are
differentiable everywhere and therefore have a tangent at any of
their points. In 2011  the problem of differentiability for
 self-affine curves with 2
generators and the problem of existence of tangent subspaces for
self-similar sets found their  exhaustive solution in the paper of
Ch.Bandt and A.Kravchenko \cite{BKra}.

In this note we study the projections of  self-similar Jordan arcs
in $\rd$ to the real line along families of parallel hyperplanes.
 Analysing the case when there is a bijective
projection of a self-similar Jordan arc $\ga$  to a straight line
segment, we show that this is possible only when the arc is a
straight line segment itself.

\begin{thm}\label{parallel}
Let $\ga$ be a self-similar Jordan arc in $\rd$. Suppose there is
such hyperplane $\sa$, that for any $x\in\ga$ the parallel copy of
$\sa$ passing through the point $x$ intersects $\ga$ only once,
then $\ga$ is a straight line segment.
\end{thm}

Really, we prove a much more general statement, in which
transverse hyperplanes $\sa(x)$ at different points $x$ of $\ga$
need not be parallel to each other, and transversality is
understood in the sense of Definition \ref{weaktransv}:

\begin{thm}\label{nonparallel}
Let $\ga$ be a self-similar Jordan arc in $\rd$. Suppose there is
a dense subset $D\IN\ga$ such that for any $x\in D$ there is  a
hyperplane $\sa$, which is weakly transverse to $\ga$ at the point
$x$, then $\ga$ is a straight line segment.
\end{thm}

The proof is based on a simple and almost obvious observation
(Theorem \ref{mz dim1}), that  the invariant set of a multizipper
of similarity dimension 1 is always a collection of straight line
segments. We prove it in Section 3.

The author is thankful to V.V.Aseev for numerous fruitful
discussions of the topic.

\section{Preliminaries.}

We give some definitions  needed in current paper. Some of them
are slightly different from generally accepted ones, but they are
best fit for our further argument.

{\bf Self-similar  arcs.} A contraction similarity $S$ in $\rd$ is
a map of the form $S(x)=q\cdot O(x-x_0)+x_0$, where $x_0$ is the
fixed point of $S$, $q\in(0,1)$ is the contraction ratio, and $O$
is the orthogonal transformation called the \emph{orthogonal part}
of $S$.

Let $\eS=\{S_1,...,S_m\}$ be a system of contraction similarities
in $\rd$. A compact set $K$ is called the \emph{invariant set} or
the \emph{attractor} of the system $\eS$, if
$K=\bigcup\limits_{i=1}^m S_i(K)$. If this invariant set is an arc
$\ga$ we call $\ga$ a \emph{self-similar arc} defined by the
system $\eS$.

We denote the semigroup generated by $S_1,...,S_m$, by $G(\eS)$.

{\bf Directed multigraphs.} A \emph{directed multigraph} (or
\emph{digraph}) $\Ga$ is defined by a set of \emph{vertices}
$V(\Ga)$, a set of \emph{edges} $E(\Ga)$ and maps $\al,
\om:E(\Ga)\to V(\Ga)$. Here $\al(e)$ is the \emph{beginning} of
the edge $e$ and $\om(e)$ is its \emph{end}.

By $E_{uv}$  we denote the set of all edges $e\in E$ for which
$\al(e)=u, \om(e)=v$, and by $E_u=\bigcup\limits_{v\in V}E_{uv}$
--- the set of all edges with the starting point at $u$.

To make the further argument more convenient, the set $V$ will be
supposed to be equal to $\{1,2,...,n\}$, where $n=\# V$. In this
case $u\in V$ means the same as $1\le u\le n$. We also denote the
numbers $\#E_{uv}$ by $m_{uv}$ and $\#E_{u}$ by $m_{u}$.

A \emph{path} $\sa$ from  a vertex $\al(e_1)=u$ to $\om(e_n)=v$ in
a digraph $\Ga$ is a sequence of edges $\sa=e_1 e_2\ldots e_n$,
with $\om(e_i)=\al(e_{i+1})$ for every $1\leq i\leq n-1$. The set
of all paths $\sa$ of the length $n$ with the beginning  $u$ and
the end  $v$ is denoted by $E^{(n)}_{uv}$ and
$E^{(*)}_{uv}=\bigcup\limits_{n=1}^\8 E^{(n)}_{uv}$ is the  set of
all paths from $u $ to $v$.

 A digraph $\Ga$ is \emph{strongly
connected} if for every two vertices $u$ and $v$ it has a path
from $u$ to $v$.

\bigskip

{\bf Graph-directed   systems of contraction similarities.}

 A \emph{graph-directed   system of contraction similarities} $\eS$ with
  \emph{structural graph} $\Ga=\langle V, E, \al, \om\rangle$ is a
  finite collection of metric spaces $\{X_v\}_{v\in V}$,
  together with a
  collection of contraction similarities $\{S_e:X_{\om(e)}\to
  X_{\al(e)}\}_{e\in E}$.

  We denote the contraction ratios of the similarities by $q_e=\Lip(S_e)$.

 Throughout  this paper all the  spaces $X_u$ will be
different copies  of the  same space $\rd$ for certain $d$.

A {graph-directed   system of similarities} $\eS$ is called
\emph{regular}, if  its structural  graph $\Ga$ is  strongly
connected.

  A
  finite collection of compact  subsets $\{K_v\}_{v\in V}$, is
  called the \emph{invariant set}, or the \emph{attractor} of the system $\eS$, if for every $v\in V$
  \begin{equation}
    K_u=\bigcup_{\al(e)=u} S_e(K_{\om(e)}).
  \end{equation}

The sets $\{K_u\}_{ u\in V}$ are called the  \emph{components of
the attractor} of the  system $\eS$.

We use the following definition of a similarity dimension of
graph-directed system of similarities \cite{Ed MTFG},\cite{MW}:

\begin{dfn}
Let $\eS$ be a regular graph-directed system of similarities with
a structure graph $\Ga=\langle V, E, \al, \om\rangle$. For each
positive real number $s$, let $\bf B(s)$ be the matrix (with rows
and columns indexed by V ) with entry ${\bf B_{uv}}(s) =
\sum\limits_{e\in E_{uv}} q_e^s$ in row $u$ column $v$. Let
$\Phi(s) = r({\bf B}(s))$  be the spectral radius of ${\bf
B}(s)$.   The unique solution $s_1\ge 0$ of $\Phi(s) = 1$ is the
similarity dimension of the system $\eS$.
\end{dfn}

\section{Multizippers of similarity dimension 1.}

A method of construction of self-similar curves, used by many
authors \cite{Serp, Levy, Hut} was studied in 2002 by V.V.Aseev
\cite{Ase} as a zipper construction. This construction proved to
be an efficient tool in the investigation of geometrical
properties of self-similar curves and continua \cite{ATK}. Its
graph-directed version was  introduced by the author in 2006 and
was called a multizipper construction; it gives a complete
description of self-similar Jordan arcs in $\rd$ \cite[Theorem
4.1]{Atet1}:

\begin{thm}
Let $\eS$ be a regular graph directed system of similarities in
$\rd$ with Jordan attractor $\vec\ga$. If one of the components
$\ga_u$ of the attractor $\vec\ga$ is different from a straight
line segment, then there is a multizipper $\eZ$ such that the set
of the components of the attractor of $\eZ$ contains each of the
arcs $\ga_u$.
\end{thm}

{\bf Definition of a multizipper.} Consider a graph-directed
system $\eZ$ of similarities  with structural graph $\Ga$, which
satisfies the following conditions:

{\bf MZ1.} In each of  the  spaces $X_u, u\in V$,   a chain of
points $\{z_0^{(u)},...,z_{m_u}^{(u)}\}$,  is specified. These
chains are   defined in such a way that
$$\|z_{i}^{(u)}-z_{i-1}^{(u)}\|<\|z_{m_v}^{(v)}-z_{0}^{(v)}\|\ \ \mbox{for any\ \ } u,v\in V, i=1,...,m_u.$$

{\bf MZ2.} There is a bijection $\epsilon $ from the  set of all
pairs $\{(u,i), u\in V, 1\le i\le m_u\}$ to the  set $E$.

{\bf MZ3.} For any pair $(u,i)$,  the map $S_e$, corresponding to
the  edge  $e=\epsilon(u,i)$ with $v=\om(e)$, sends two-point set
$\{z_{0}^{(v)},z_{m_v}^{(v)}\}$ to the set
$\{z_{i-1}^{(u)},z_{i}^{(u)}\}$.

The graph-directed system $\eZ$, satisfying the conditions {\bf
MZ1--- MZ3} is called a \emph{multizipper with structural graph
$\Ga$ and node points $z^{(u)}_i.$}

\vspace{.7cm}

Let $L^{(u)}$ be the polygonal line specified by the sequence
$\{z_0^{(u)}, z_1^{(u)},...,z_{m_u}^{(u)}\}$ of the nodes  of the
multizipper $\eZ$. Denote the distance
$||z_{m_u}^{(u)}-z_{0}^{(u)}||$ by $l_u$. Observe that if
$S_e(\{z_{0}^{(v)},z_{m_v}^{(v)}\})=\{z_{i-1}^{(u)},z_{i}^{(u)}\}$,
then  $||z_{i-1}^{(u)} - z_{i}^{(u)}||=q_el_v$. So, the length of
the polygonal line $L^{(u)}$ is  equal to
$\sum\limits_{v=1}^n\sum\limits_{e\in E_{uv}}q_{e} l_v.$

\vspace{.7cm}


\begin{thm}\label{mz dim1}
Let $\eZ$ be a regular self-similar multizipper whose similarity
dimension is 1. Then all the  components $\ga^{(u)}$ of its
invariant set are line  segments.
\end{thm}

\dok. Suppose there is a component $\ga^{(u)}$ of the attractor of
$\eZ$, which is not a line segment. Since $\eZ$ is regular, for
any $v\in V$ there is a path $\sa=e_1...e_k\in E^{(*)}_{vu}$, so
the similarity $S_\sa=S_{e_1}\cdot...\cdot S_{e_k}$ maps the arc
$\ga^{(u)}$ to a subarc of $\ga^{(v)}$. Therefore, each
$\ga^{(v)}$ is also different from a straight line.

Then, choosing appropriate refinement of the multizipper   $\eZ$,
we may suppose that all the polygonal lines $L^{(u)}$ are
different from a straight line. For  each component $\ga^{(u)}$ we
have:
$$\ga^{(u)}=\bigcup\limits_{v=1}^n\bigcup\limits_{e\in E_{uv}}S_{e}(\ga^{(v)}).$$

The similarity dimension of the multizipper $\eZ$ is equal to such
value of a parameter $s$, that the  spectral radius of  the matrix
$\bf B(s)$ whose  entries are $B_{uv}(s)=\sum\limits_{e\in
E_{uv}}q_{e}^s$, is equal to 1.

So, the spectral  radius of the matrix $\bf B(1)$ with entries
$B_{uv}(1)=\sum\limits_{e\in E_{uv}}q_{e}$ is equal to 1.

Since all the polygonal lines $L^{(u)}$ are not straight lines,
they obey the inequality
$$l_u<\sum\limits_{v=1}^n\sum\limits_{e\in E_{uv}}q_{e} l_v=({\bf B(1)} \vec l)_u. $$

Therefore, for a vector $\vec l=(l_1,...,l_n)$ and for the matrix
$\bf B(1)$ we have the inequality $$\mathop{min}\limits_{1\le u\le
n}\dfrac{({\bf B} \vec l)_u}{l_u}>1.$$

The structural graph of the system $\eZ$ is  strongly connected.
Then the matrix $\bf B(1)$ is a positive irreducible matrix.
According to \cite[Remark 4, \S2, Ch.XIII]{Gan} its spectral
radius is equal to
$$r=\mathop{max}\limits_{\vec l\neq 0}\mathop{min}\limits_{1\le u\le n}\dfrac{({\bf B}\vec l)_u}{l_u}.$$

So, if $r=1$, then for any $\vec l$, there is such $u$, that
$\dfrac{({\bf B}\vec l)_u}{l_u}\le 1$.

The contradiction shows that all $L^{(u)}$ are straight line
segments, so all $\ga^{(u)}$ are  straight line segments
too.$\blacksquare$

\section{Theorem on transverse hyperplanes.}

{\bf Jordan arcs and transverse hyperplanes.}
 Let $\ga:[0,1]\to \rr^d$ be a Jordan arc in $\rr^d$. For any point
  $x=\ga(t)$ we  define the half-open subarcs  $\ga^+_x=\ga((t,1])$ and  $
 \ga^-_x=\ga([0,t))$.

 Let $x,y\in\ga, $ and $y\in \bar \ga^+_x$.  We denote  the open
subarc $\ga^+_x\cap\ga^-_y$ by $(x,y)$ and
$\bar\ga^+_x\cap\bar\ga^-_y$
 by $[x,y]$.

A hyperplane containing  the origin 0 is denoted by $\sa$, while
$V^+(\sa)$ and $V^-(\sa)$ are open half-spaces, defined  by $\sa$.
A hyperplane parallel to $\sa$ and  containing $x$ is denoted by
$\sa(x)$ or $\sa+x$. The open half-spaces defined by $\sa(x)$ are
denoted by $V^+(\sa,x)$ and $V^-(\sa,x)$ or $V^+(\sa)+x$ and
$V^-(\sa)+x$.

 \begin{dfn}\label{weaktransv}
We say a hyperplane $\sa$  is   \emph{weakly transverse} to the
arc $\ga$ at  the point $x$, if $\ga_{x}^{+}\IN \bar V^+(\sa,x)$,
$\ga_{x}^-\IN \bar V^-(\sa,x)$.

We say a hyperplane $\sa$  is   \emph{ transverse} to the arc
$\ga$ at  the point $x$, if $\ga_{x}^{+}\IN   V^+(\sa,x)$,
$\ga_{x}^-\IN  V^-(\sa,x)$.
\end{dfn}

\vspace{.5cm}

{\bf The cones $Q^+$ and $Q^-$.} By $Q^+(x,y)$ \ (respectively,
$Q^-(x,y)$)\ we denote the intersection of
 all closed half-spaces
 $\bar V^+(\sa,z)$\  ( resp.
$\bar V^-(\sa,z)$)\  corresponding to the hyperplanes $\sa(z)$,
weakly transverse to $\ga$ at  the points $z\in[x,y]$. These sets
 are  convex and  closed and  they satisfy the relations $$\ga^+(y)\IN
Q^+(x,y)\ \ \mbox{ and } \ga^-(x)\IN Q^-(x,y).$$

 Taking $x=y$ we come to the sets $Q^+(x)$ \ (  $Q^-(x)$) \  which are the intersections
  of all closed half-spaces $\bar V^+(\sa,x)$ \ ($\bar
  V^-(\sa,x)$) \
corresponding to hyperplanes $\sa(x)$, weakly transverse to $\ga$
at the point $x$. We can also consider the  set $Q^+(x)\cup
Q^-(x)$ as the intersection of all unions $Q^+_i\cup Q^-_i $ of
pairs of convex closed cones symmetric with respect  to $x$ which
satisfy relations $\ga^+(x)\IN Q^+_i$ and $\ga^-(x)\IN Q^-_i$.

\begin{lem}\label{limsa}
Let $\ga$ be a Jordan arc in $\rr^n$. Suppose a sequence of points
$x_n\in\ga$ converges to a point $x_0$, while a sequence of
hyperplanes $\sa_n$, weakly transverse to  $\ga$ at points $x_n$,
converges to a hyperplane $\sa_0$. Then $\sa_0$ is weakly
transverse to  $\ga $ at the point $x_0$.

\end{lem}

\dok. For  any $n$,  $\bar \ga^+_{x_n}\IN V^+(\sa_n,x_n)$. Since
$x_n\to x_0$, $\sa_n$ converge to $\sa_0$ if and only if
$\sa_n(x_n)$ converge to $\sa_0(x_0)$. Taking the closed
half-spaces, corresponding to $\sa_n(x_n)$, we get
$\lim\limits_{n\to\8}\bar V^+(\sa_n,x_n)=\bar V^+(\sa_0,x_0)$. At
the same time, $\lim\limits_{n\to\8}\bar \ga^+_{x_n}=\bar
\ga^+_{x_0}$. Therefore,
 $\bar\ga^+_{x_0}\IN \bar V^+(\sa_0,x_0)$. The  same  way we
get   $\bar\ga^-_{x_0}\IN \bar V^-(\sa_0,x_0)$. \vse

\medskip

Denote by $\Sa(x)$ the set of all hyperplanes, weakly transverse
to the arc $\ga$ at the point $x\in\ga$. This set is a compact
subset of $\bf{RP}^d$. It follows  from the Lemma \ref{limsa},
that
$$\Sa(x)\NI \mathop{\lim\sup}\limits_{y\to x,y\in\ga}\Sa(y).$$

This inclusion implies that the cones $Q^+(x)$  and  $Q^-(x)$
satisfy the following semicontinuity condition:

\begin{lem}\label{liminf}
Let $\ga$ be a Jordan arc in  $\rr^d$ and $x\in\ga$. Then,
 $$Q^+(x)\IN\mathop{\lim\inf}\limits_{y\to x,y\in\ga}Q^+(y). \eqno$$

\end{lem}

 \dok.
Since $$Q^+(x)=\bigcap\limits_{\sa\in\Sa(x)}\bar V^+(\sa,x),$$
 using basic properties of upper and lower limits\cite[\S 29]{Kur}, we
can write $$Q^+(x)=(\bigcup\limits_{\sa\in\Sa(x)}
V^-(\sa,x))^\mathsf{c}\IN\ \ (\mathop{\lim\sup}\limits_{y\to
x,y\in\ga}\bigcup\limits_{\sa\in\Sa(y)} V^-(\sa,y))^\mathsf{c}=$$
 $$=\mathop{\lim\inf}\limits_{y\to
x,y\in\ga}(\bigcup\limits_{\sa\in\Sa(y)} V^-(\sa,y))^\mathsf{c}
=\mathop{\lim\inf}\limits_{y\to x,y\in\ga}Q^+(y).\blacksquare$$

\begin{lem}\label{weak to strong}
Let $\ga$ be a self-similar Jordan arc. If for any $x\in\ga$ there
is a hyperplane, weakly transverse to $\ga$ at the point $x$, then
there is a  hyperplane $\sa$, which is transverse to $\ga$ at any
point $x\in\ga$.
\end{lem}

\dok. Suppose the affine hull of $\ga$ is $\rd$ so it is not
contained in a hyperplane.

Take some $\da>0$.

Consider the  family of all the cones $A=\{Q^+(x),x\in\ga\}$.
Taking the parallel copy of each cone $Q^+(x)$ having the vertex
at the center $0$ of the unit ball $B\IN\rr^d$, we denote its
intersection with the ball $B$ by $Q(x)$. This turns the family
$A$ to a   subset of the hyperspace $Conv(B)$ of compact convex
subsets of the unit ball $B$. Observe that the inclusion
$Q(x)\IN\mathop{\lim\inf}\limits_{y\to x,y\in\ga}Q(y)$ in the
statement of Lemma \ref{liminf} holds for the cones $Q(x)$ as
well.

Let $S $ be a contraction similarity, for which $S(\ga)\IN\ga$.
Let $x_0$ be its fixed point. Let $O $ be the orthogonal part of
the similarity $S $.

Since $Q(x_0)\IN \mathop{\lim\inf}\limits_{x\to x_0}Q(x)$, there
is an open subarc $(y,z)\ni x_0$ such that for any $x\in(y,z)$,
the cone
  $ Q(x_0)$ is contained in $\da-$neighborhood  $N_\da(Q(x))$ of a cone $Q(x)$.

For some sufficiently large $k$, the subarc $S_i^k(\ga)$ is
contained in $ (y,z)$. Then for any $\xi\in\ga$, the point
$x=S_i^k(\xi)$ lies in $ (y,z)$ and $N_\da(Q(x))\NI Q(x_0)$. Since
$Q(x)=O_i^k(Q(\xi))$, and $Q(x_0)=O_i^k (Q(x_0))$ and  $O_i^k$ is
an isometry, $N_\da(Q(\xi)) $ must also contain $ Q(x_0)$.

Thus, if  $S:\ga\to\ga$ is a similarity and  $\fix(S)=x$, then for
any $\da>0$ and any $\xi\in\ga$, $N_\da(Q(\xi))\NI Q(x)$.
Therefore, $Q(\xi)\NI Q(x)$ for all $\xi\in\ga$. If we take for
$\xi$ a fixed point of  some other similarity
 $S':\ga\to\ga$, we get that
$Q(\xi)=Q(x)$. Thus, the minimal cone $Q(x)$ is the same, no
matter which fixed point $x$ we choose, and we denote it by $Q$.
If $x$ is not a fixed point of any $S\in G(\eS)$, then $Q(x)\IN
Q$. If $\sa(x)$ is a support hyperplane to the cone $Q^+(x)$ at
some fixed point $x$, then for any $\xi\in\ga$ parallel hyperplane
$\sa(\xi)$ is a support hyperplane for $Q^+(\xi)$ and is thus
weakly transverse to  $\ga$ at the point $\xi$.

Suppose for some $x$ and $w$ in $\ga$, $ w\in\dot\ga^+(x)$ and
$w\in\sa(x)$. Then $\sa(x)=\sa(w)$ and $V^+(\sa,x)=V^+(\sa,w)$. By
weak transversality of $\sa$ at the points $x$ and $w$,  the
subarc $[x,w]$ lies in $V^+(\sa,x)\cap V^-(\sa,w)=\sa$. Then the
whole arc $\ga$ lies in a hyperplane. The contradiction shows that
the hyperplanes parallel to $\sa(x)$ are transverse to $\ga $ at
any point.\vse

\begin{lem}\label{perplane}
Let $\ga$ be a self-similar Jordan arc, which has a hyperplane
 transverse to $\ga$ at any of its points. Then
there is such transverse hyperplane $\sa$, that for any similarity
 $S_i\in \eS$, $O_i(\sa)=\sa$.
\end{lem}

\dok.   Let $G_O$ be a group generated by orthogonal parts $O_i$
of the similarities $S_i\in \eS$. For any $O\in G_O$, the image
$O(Q)$ is either $ Q$  or $-Q$. The space $\rr^d$ is  a direct sum
of two orthogonal subspaces
  $X_0\oplus X_1$, where  $X_0$ is the space of all such $x$ that for any $O\in
  G_O$,
    $O(\{x,-x\})=\{x,-x\}$ and $X_1=X_0^\bot$.

  Consider the intersection $X_0\cap  Q$. This intersection is a
  convex cone $Q'$ in $X_0$. Take a support hyperplane $Y$ to the cone $Q'$
  at the point $0$ in the space $X_0$. Then $Y+X_1$  is a support hyperplane for $Q$ in $\rd$.

  Suppose contrary. Then there is some  $z\in (Y+X_1)\cap \dot Q$. The point $z$
  has unique representation in the  form $z=x+y$, where $x\in X_1, \ x\neq 0$ and $y\in Y$.
Consider the convex hull $W$ of the orbit $G_O(x)$.  It's
barycenter is fixed  by the group $G_O$, therefore it is $0$. Then
the barycenter of the convex hull of the orbit $G_O(z)$ is $y$.

 Take a ball $B(z,\e)\IN Q$. The convex hull of a set
  $\bigcup\limits_{O\in G_O}O(B(z,\e))$ contains the ball $B(y, \e)$, therefore $y\in\dot
  Q$, which is impossible. So $\dot Q\cap(Y+X_1)=\0$.

At the same time, for any $O\in G_O$ the transformation $O$ sends
the hyperplane $Y+X_1$ to itself.\vse

\vspace{1cm}

{\bf The proof of Theorem \ref{nonparallel}.}

 Let $\ga$ be a self-similar
Jordan arc, which is not a line segment. By Theorem 4.1 in
\cite{Atet1}, the arc $\ga$ may be represented as a component
$\ga^{(u)}$ of the invariant set of some multizipper $\eZ$ , for
which the maps $S_e, e\in E$ are the elements of the semigroup
$G(\eS)$. Let $z^{(u)}_i$ be the node points  and $\Ga=\langle
V,E,\al,\om \rangle$ be the structural graph of $\eZ$. Passing, if
necessary, to a subarc of $\ga$, we may suppose that the graph
$\Ga$ is  strongly connected and the multizipper $\eZ$ is regular.

If $\ga$ contains such dense subset $D\IN \ga$, that for any $x\in
D$ there is a hyperplane $\sa(x)$, weakly transverse to $\ga$,
then by Lemma \ref{limsa}, such hyperplane $\sa(x)$ exists for any
$x\in\ga$. By Lemma \ref{weak to strong}, there is a hyperplane
$\sa$, transverse to $\ga$ at any $x\in\ga$.

 By  Lemma \ref{perplane} there is a hyperplane
$\sa$, transverse to $\ga$ at  any of its points, which is
preserved by any of $O_i\in G_O$. Then the duplicates of $\sa$ are
transverse to the  components $\ga^{(u)}, u\in V$ of the attractor
 of the multizipper $\eZ$ at any of their points and are preserved
by the orthogonal parts $O_e$ of the similarities $S_e$.

Let $\La^{(u)}$ be a line, orthogonal to $\sa$ in the copy
$X^{(u)}$ of the space $\rd$. Let $\ga^{(u)}$ be the component of
the invariant set of $\eZ$ lying in $X^{(u)}$. Consider the
orthogonal projection $\pi$ of each arc $\ga^{(u)}$ to the
$\La^{(u)}$.

Since the similarities $S_e$ send the hyperplanes, parallel to
$\sa$, to the hyperplanes, parallel to $\sa$, for each similarity
$S_e\in\eZ$, $S_e:\ga^{(v)}\to \ga^{(u)}$ there is a similarity
$\hat S_e:\La^{(v)}\to \La^{(u)}$, satisfying the condition
$$\pi\circ S_e=\hat S_e\circ\pi. \ \ \ \eqno \label{equivar}$$
Due to this  condition each map $\hat S_e$ sends the set
$\{\pi(z_{0}^{(v)}),\pi(z_{m_v}^{(v)})\}$ to the set
$\{\pi(z_{i-1}^{(u)}),\pi(z_{i}^{(u)})\}$.

The system $\hat\eZ$ is a linear multizipper with node points
$\hat z^{(u)}_i=\pi(z_{i}^{(u)})$.

Since for  any $S_e$, $\Lip(\hat S_e)=\Lip(S_e)$ the similarity
dimension of the multizipper $\eZ$  is  equal to the similarity
dimension of $\hat\eZ$ and therefore it is  equal to 1. By Theorem
\ref{mz dim1}, its invariant set is a collection of straight line
segments. $\blacksquare$


\begin{thebibliography}{99}

\bibitem{Ase}
V.~V.~Aseev,  On the regularity of self-similar zippers. --
  "The 6-th Russian-Korean Int. Symp. on Sci. and Technology.
  KORUS-2002 (June 24-30, 2002. Novosibirsk State Techn. Univ., Russia).
  Part 3 (Abstracts)", p.167

\bibitem{ATK}
V.~V.~Aseev, A.~V.~Tetenov and A.~S.~Kravchenko, On Selfsimilar
Jordan Curves on the Plane., Siberian Math.J. 44, No. 3 (2003),
379-386.




\bibitem{BKra} Ch.~Bandt, A.~Kravchenko. Differentiability of fractal curves,  Nonlinearity, 2011, 24, pp.~2717-
2728.

 \bibitem{Grig}  L.~Bartholdi, R.~I.~Grigorchuk,
 V.~V.~Nekrashevych.
From fractal groups to fractal sets., Arxiv.org preprint
math.GR/0202001,2002.



\bibitem{Ed MTFG} G.~A.~Edgar,  Measure, Topology, and Fractal Geometry, Springer-Verlag, 1990.

 \bibitem{EdgDas}  G.~A.~Edgar, M.~Das.
Separation properties  for  graph-directed self-similar fractals.,
Top.~appl.~,2005, V.152, p.138--156.


 \bibitem{Fal}
 K.~J.~Falconer. Fractal geometry: mathematical foundations and
applications . --  J.~Wiley and Sons, New York, 1990.

\bibitem{Gan}
 F.~R.~Gantmacher, The Theory of Matrices., Volume 2, Chelsea, 1959,
Chapter XIII: Matrices with non-negative elements.




\bibitem{Hut}
J.~Hutchinson. Fractals and self-similarity.   Indiana Univ. Math.
J., 1981, V.~30, No.~5,  pp.~713--747.





\bibitem{Koch}
H.~von~Koch,  Sur une courbe continue sans tangente, obtenue par
une construction geometrique elementaire.// Archiv for Matemat.,
Astron. och Fys., 1904, V.~1, p.~681--702.

\bibitem{Krav} A. Kravchenko, Smooth self-affine curves (in Russian), Preprint
No. 161, Sobolev Institute of Mathematics, Siberian Branch of the
Russian Academy of Sciences, 2005.

\bibitem{Kur} K.~Kuratowski. Topology. Volume I.   Academic Press/Polish Scientific
Publishers, New York/London/Warszawa 1966.


\bibitem{Levy}
 P.~Levy. Les courbes planes ou gauches et les surfaces composees de parties
semblables au tout., J.~Ecole Polytechn., III. Ser. 1938. V.~144,
P.~227--247 et 249--291


 \bibitem{MW}   R.~D.~Mauldin, S.~C.~Williams. Hausdorff dimension in graph
directed constructions., Trans. Amer. Math. Soc. 1988, V.~309,
pp.~811--829.

\bibitem{Serp} W.~Sierpinski.  Sur une courbe dont tout point est un point de
ramification//
 Compt. Rendus Acad. Sci. Paris, 1915, V.~160, pp.~302--305.


\bibitem{Atet1}  A.~V.~Tetenov. Self-similar Jordan arcs and
 graph-directed systems of similarities, Siberian.Math.J., 47, No.5 (2006),
 pp.~940-949.














\end{thebibliography}
\end{document}